\documentclass[final,nomarks]{dmtcs-episciences}
\usepackage[utf8]{inputenc}
\usepackage{subfigure}
\usepackage{graphicx}
\usepackage[english]{babel}
\usepackage{amssymb}
\usepackage{amsmath}
\usepackage{latexsym}
\usepackage{amsthm}
\usepackage{tikz} 
\usepackage{hyperref}
\usetikzlibrary{cd}
\usetikzlibrary{decorations.pathreplacing}
 \theoremstyle{plain}
 \newtheorem{thm}{Theorem}[section]
 \newtheorem{cor}[thm]{Corollary}
 \newtheorem{lem}[thm]{Lemma}
 \newtheorem{prop}[thm]{Proposition}

 \theoremstyle{definition}
 
 \newtheorem{example}[thm]{Example}

 \theoremstyle{remark}
 
\author{Marilena Barnabei 
\and
Flavio Bonetti \\
\and
Niccol\`o Castronuovo 
\and
Matteo Silimbani}

\title[Consecutive patterns]{Consecutive patterns in restricted permutations and involutions  \thanks{This work was partially supported by the University of Bologna, funds for selected research topics.}}
\affiliation{
  Dipartimento di Matematica, Universit\`a di Bologna, Bologna, 40126, ITALY}
\received{2019-2-7}
\revised{2019-5-17}
\accepted{2019-5-17}

\begin{document}
\publicationdetails{21}{2019}{3}{21}{5175}
\maketitle
\begin{abstract}
It is well-known that the set $\mathbf I_n$ of involutions of the symmetric group $\mathbf S_n$ corresponds bijectively  - by the Foata map $F$  - to the set of $n$-permutations that avoid the two vincular patterns $1\underline{23}$ and $1\underline{32}.$ We consider a bijection $\Gamma$ from the set $\mathbf S_n$ to the set of histoires de Laguerre, namely, bicolored Motzkin paths with labelled steps, and study its properties when restricted to $\mathbf S_n(1\underline{23},1\underline{32}).$ In particular, we show that the set $\mathbf S_n(\underline{123},{132})$
of permutations that avoids the consecutive pattern $\underline{123}$ and the classical pattern $132$ corresponds via $\Gamma$ to the set of Motzkin paths, while its image under $F$ is the set of restricted involutions $\mathbf I_n(3412).$
We exploit these results  to determine the joint distribution of the statistics des and inv over 
 $\mathbf S_n(\underline{123},{132})$ and over $\mathbf I_n(3412).$
 Moreover, we determine the distribution in these two sets of every consecutive pattern of length three. To this aim, we use a modified version of the well-known Goulden-Jacson cluster method.
\end{abstract}

\section{Introduction}

In 1994 De Medicis and Viennot \cite{DeM} introduced the definition of histoire de Laguerre, namely, a pair $(d,l),$ where $d$ is a Motzkin path of length $n$
whose horizontal steps may have two different colors  and $l=(l_1,\ldots ,l_n)$ is a sequence of non-negative integers with suitable constraints. 

Many bijections are present in the literature between the  set $\mathcal H_n$ of histoires de Laguerre and the symmetric group $\mathbf S_n$ 
(see \cite{Rand} for a survey on this topic),  as well as between a specific subset $\mathcal L_n$  of
$\mathcal H_n$ and the set $\mathbf I_n$ of involutions in $\mathbf S_n$  (see, \textit{e.g.}, \cite{Ba6}). 
 More recently,  Claesson \cite{Claesson} proved that the set $\mathbf I_n$ corresponds bijectively - via the classical Foata map - with the set $\mathbf S_n(1\underline{23},1\underline{32})$ of $n$-permutations  that avoid the  two vincular patterns $1\underline{23}$ and $1\underline{32}.$ 
 
In the first part of the  present paper we consider  the bijection  $\Gamma$ between the set of $n$-permutations and the set of historires de Laguerre described in \cite{BBCS}, which is essentially   the bijection defined in \cite[p. 256]{CSZ}, up to the reverse. This last bijection is in turn a slightly modified version of the well-known Fran\c{c}on-Viennot bijection \cite{Francon}. We exploit $\Gamma$ 
to connect the sets $\mathbf S_n(1\underline{23},1\underline{32}),$ $\mathbf I_n$ and $\mathcal L_n.$ More precisely we prove that the image under $\Gamma$ of the set $\mathbf S_n(1\underline{23},1\underline{32})$ is precisely the set $\mathcal L_n,$ namely, the set  consisting of (uncolored) Motzkin paths whose down steps are labelled with an integer that does not exceed their height.
Furthermore, we show that the bijection $\Psi$ defined in \cite{Ba6} is nothing but the composition of the map $F$ with $\Gamma.$ Finally, the set of permutations avoiding the consecutive pattern $\underline{123}$ and the classical pattern $132$ is mapped by $\Gamma$  onto the   set of unlabelled Motzkin paths and is mapped by $F^{-1}$
onto the set of involutions avoiding the classical pattern $3412.$ 

In the second part of the paper we exploit the properties of the maps $\Gamma$ and $\Psi$  to study in parallel some statistics over the two sets $\mathbf S_n(\underline{123},132)$
and $\mathbf I_n(3412).$ In particular, in both cases we determine the joint distribution of inversions and descents, as well as the distribution of the occurrences of every consecutive pattern of length three.

In many situations we take advantage of a particular instance of the Goulden-Jackson cluster method \cite{Goulden} for Motzkin paths. For the sake of completeness we describe this method in the Appendix.

\section{The bijections}\label{sectone}

A \emph{Motzkin path} of length  $n$ is a lattice path 
starting in $(0,0),$ ending in $(n,0),$ consisting of \emph{up steps}  $U$ of the form $(1,1),$  \emph{down steps} $D$ of the form $(1,-1)$ and \emph{horizontal steps} $H$ of the form $(1,0)$
and lying weakly above
the $x$-axis.

As usual, a Motzkin path can be identified with a \emph{Motzkin word}, namely, a word $w=d_1d_2\ldots d_{n}$ of length $n$ in the alphabet $\{U,D,H\}$ with the constraint that
the number of occurrences of the letter $U$ is equal to the number of occurrences of the letter $D$ and, for every $i,$
the number of occurrences of $U$ in the subword $d_1d_2\ldots d_i$ is not smaller than the number of occurrences of $D.$
In the following we will not distinguish between a Motzkin path and the corresponding word.

Now we consider the set of \emph{bicolored } Motzkin paths, defined as Motzkin paths whose horizontal steps have two possible colors $c_1$ and $c_2,$ 
such that horizontal steps lying on the $x$-axis cannot be colored with the color $c_2.$ We will denote by $H$ a horizontal step colored by $c_1$ and with $\widetilde H$ a horizontal step colored by $c_2.$
 It is well known that bicolored Motzkin paths are counted by Catalan numbers (see \cite{St3}).

We will denote by $\mathcal M_n$ and $\mathcal {CM}_n$ the sets of Motzkin paths of length $n$ and bicolored Motzkin paths of length $n,$  respectively.

 We associate to every $d = d_1d_2\ldots d_{n} \in \mathcal {CM}_n$ the $n$-tuple $h(d)=(h_1,h_2,\ldots,h_n),$ where 
for every $i=1,\ldots ,n,$ the integer $h_i$ is defined as 
\begin{displaymath}\begin{cases}
\mbox{the y-coordinate of the ending point of the step } d_i & \mbox{ if $d_i$ is $D$} \\
\mbox{the y-coordinate of the starting point of the step } d_i & \mbox{otherwise}
\end{cases}\end{displaymath}

  We will call the integer $h_i$ the \emph{height} of the step $d_i,$ and the $n$-tuple  $h(d)$ the \emph{height list} of $d.$

\begin{example}
Consider the bicolored Motzkin path $d=UUD\widetilde{H}DH,$ namely,
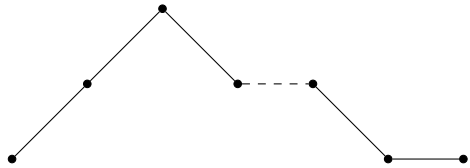
\begin{displaymath}
\begin{tikzpicture}%[every node/.style={draw,shape=circle,minimum size=1mm, inner sep=0mm, fill=black}]
\node[draw,shape=circle,minimum size=1mm, inner sep=0mm, fill=black] (6) at (6,0) {};
   \node[draw,shape=circle,minimum size=1mm, inner sep=0mm, fill=black] (5) at (5,0) {};
   \node[draw,shape=circle,minimum size=1mm, inner sep=0mm, fill=black] (4) at (4,1) {};
   \node[draw,shape=circle,minimum size=1mm, inner sep=0mm, fill=black] (3) at (3,1) {};
   \node[draw,shape=circle,minimum size=1mm, inner sep=0mm, fill=black] (2) at (2,2) {};
   \node[draw,shape=circle,minimum size=1mm, inner sep=0mm, fill=black] (1) at (1,1) {};
   \node[draw,shape=circle,minimum size=1mm, inner sep=0mm, fill=black] (0) at (0,0) {};
   \draw[] (0) -- (1) -- (2) -- (3)    (4) -- (5) -- (6);
   \draw[dashed] (3) -- (4);
   
    \end{tikzpicture}
\end{displaymath}
where the horizontal step with color $c_2$ is represented by a dashed line.
Then $h(d)=(0,1,1,1,0,0).$
\end{example}

\bigskip

We now describe a map from the set of permutations of length $n$ to the set $\mathcal {CM}_n$. This map is a slight modification of the map described in \cite{BBCS} in terms of valued Dyck paths.  

 Let $\pi=\pi_1\pi_2\ldots \pi_n$ be a permutation in $\mathbf S_n$ written in one-line notation.
An \emph{ascending run} in $\pi$ is a maximal increasing subword of $\pi.$ For example, the ascending runs of $346512$ are $w_1=346,$ $w_2=5$ and $w_3=12.$
Write $\pi$ as \begin{displaymath} \pi=w_1w_2\ldots w_k,\end{displaymath} where the ${w_i}'s$ are the ascending runs in $\pi.$
The first and the last  element of an ascending run of length at least two are called a \emph{head} and a \emph{tail}, respectively.
The only element of an ascending run of length one is called a \emph{head-tail}.
Every other element is called a \emph{boarder}.

Now we associate to $\pi$ a bicolored Motzkin path $d$ of length $n$ defined as follows.
For $i=1,\ldots ,n,$
\begin{itemize}
\item if $i$ is a head-tail, set $d_i=H$;
\item if $i$ is a head, set $d_i=U$;
\item if $i$ is a tail, set $d_i=D$;
\item if $i$ is a boarder, set $d_i=\widetilde{H}$.
\end{itemize}

Then $d=d_1d_2\ldots d_n.$

Obviously the correspondence $\gamma: \pi \rightarrow d$ is far from being injective. For example, both the permutations $3124$ and $1243$ in $\mathbf S_4$ correspond to the
bicolored Motzkin  path $U \widetilde{H} H D.$
In order to get a bijection, we associate to the permutation $\pi$  a pair $(d,l),$ where $d$ is the bicolored Motzkin path defined above and $l=(l_1,l_2,\ldots ,l_n)$ is the sequence of non-negative integers
\begin{displaymath} l_i=|\{j | s_j < i < t_j,\;  t_j\mbox{ precedes }i\mbox{ in }\pi\}|\end{displaymath}
where $s_j$ and $t_j$ are the first and the last element of the $j$-th ascending run of $\pi$.

We denote by $\Gamma(\pi)$ the pair $(d,l)$ associated with the permutation $\pi.$

\begin{example}\label{firstex}
Consider the permutation $\pi=826913547.$ The ascending runs of $\pi$ are $w_1=8,$ $w_2=269,$  $w_3=135$ and $w_4=47.$ 

We have $\Gamma(\pi)=(d,l),$  where $d=UU\widetilde{H}UD\widetilde{H}DHD=$
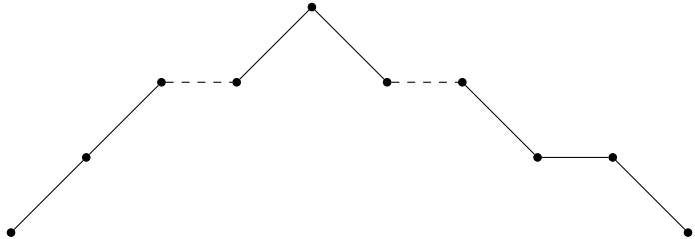
\begin{displaymath}
\begin{tikzpicture}%[every node/.style={draw,shape=circle,minimum size=1mm, inner sep=0mm, fill=black}]
 \node[draw,shape=circle,minimum size=1mm, inner sep=0mm, fill=black] (9) at (9,0) {};
 \node[draw,shape=circle,minimum size=1mm, inner sep=0mm, fill=black] (8) at (8,1) {};
   \node[draw,shape=circle,minimum size=1mm, inner sep=0mm, fill=black] (7) at (7,1) {};
   \node[draw,shape=circle,minimum size=1mm, inner sep=0mm, fill=black] (6) at (6,2) {};
   \node[draw,shape=circle,minimum size=1mm, inner sep=0mm, fill=black] (5) at (5,2) {};
   \node[draw,shape=circle,minimum size=1mm, inner sep=0mm, fill=black] (4) at (4,3) {};
   \node[draw,shape=circle,minimum size=1mm, inner sep=0mm, fill=black] (3) at (3,2) {};
   \node[draw,shape=circle,minimum size=1mm, inner sep=0mm, fill=black] (2) at (2,2) {};
   \node[draw,shape=circle,minimum size=1mm, inner sep=0mm, fill=black] (1) at (1,1) {};
   \node[draw,shape=circle,minimum size=1mm, inner sep=0mm, fill=black] (0) at (0,0) {};
   \draw[] (0) -- (1) -- (2)     (3) --  (4) -- (5)     (6) -- (7) -- (8) -- (9);
   \draw[dashed] (2) -- (3)   (5) -- (6);
   \end{tikzpicture}
\end{displaymath}
and $l=(0,0,1,2,1,0,1,0,0).$
\end{example}

We observe that the above bijection  is essentially the map defined  in \cite[p. 256]{CSZ}, up to the reverse. 

We recall that a pair $(d,l),$ where $d$ is a bicolored Motzkin path of length $n$ and $l=(l_1,\ldots ,l_n)$ is a sequence of non-negative integers, is called a \emph{histoire de Laguerre} provided that $l_i\leq h_i$ for all $1\leq i\leq n,$  where $h_i$ is the $i$-th element of the height list of $d$ (see \cite{DeM}).  We denote by $\mathcal H_n$ the set of histoires de Laguerre of length $n$.
 %Our next goal is to prove that the map $\Gamma$ is actually a bijection between $\mathbf S_n$ and $\mathcal{H}_n.$

\begin{thm}
The map $\Gamma$ is a bijection between $\mathbf S_n$ and $\mathcal{H}_n.$
\end{thm}
\proof
See \cite[Theorem 2.6]{BBCS}.
\endproof

We now describe the connection between the map $\Gamma$ defined above and a bijection $\Psi$ between the set $\mathbf I_n$ of involutions of length $n$ and labeled Motzkin paths studied in \cite{Ba6}. In order to do this, we exploit a result proved by Claesson \cite{Claesson}, namely, the fact that the classical Foata map induces a bijection between the set of involutions of length $n$ and the set of permutations of the same length that avoid two vincular patterns.

Let $\pi\in \mathbf S_n$ and $\tau\in \mathbf S_m.$ We say that $\pi=\pi_1\ldots \pi_n$ \emph{contains the pattern} $\tau=\tau_1\ldots \tau_m$ \emph{in the classical sense} if there exists an index subsequence $1\leq i_1<i_2<\ldots <i_m\leq n$ such that the words $\pi_{i_1}\pi_{i_2}\ldots \pi_{i_m}$ and $\tau_1\tau_2\ldots \tau_m$ are order isomorphic.
Otherwise, $\pi$ \emph{avoids the pattern} $\tau.$  

A \emph{vincular pattern}  is a permutation $\tau$  in $\mathbf S_m$ some of whose consecutive letters may be underlined. If
$\tau$ contains $\underline{\tau_i\tau_{i+1}\ldots \tau_j}$ as a subword then the letters corresponding to $\tau_i,\tau_{i+1},\ldots,\tau_j$  in an occurrence of $\tau$ in a permutation $\sigma$ must be adjacent, whereas there is no adjacency condition for non underlined consecutive letters (see  \cite[p. 10]{Ki}).

For example, the permutation 431256 contains two occurrences of the vincular pattern $2\underline{13}$, namely, 425 and 325.
 Note that a vincular pattern without underlined letters is a pattern in the classical sense.
On the other hand, the occurrences of a vincular pattern  all of whose letters are underlined must be formed by adjacent letters. 
The set of permutations of length $n$ that avoid the vincular pattern $\tau$ is denoted by $\mathbf S_n(\tau)$.

We now recall Claesson's result. Let $\pi$ be an involution. Write $\pi$ in \emph{standard cycle notation,} \textit{i.e.}, so that each cycle is written with its least element first
and the cycles are written in decreasing order of their least element. Define $F(\pi)$ to be the permutation obtained from $\pi$ by erasing the parentheses separating the cycles.
As an example consider $\pi =47318625\in \mathbf I_8.$ The cycle notation for $\pi$ is $(6)(5,8)(3)(2,7)(1,4) $ and $F(\pi)=65832714.$ 

In \cite{Claesson} Claesson proved that the map $F$ is a bijection between $\mathbf I_n$ and $\mathbf S_n(1\underline{32},1\underline{23})$.
It is easily seen that this last set coincides with $\mathbf S_n(\underline{132},\underline{123})$  (see \cite{ElizaldeNoy}).
On the other hand, in \cite{Ba6} the authors define a bijection $\Phi$ between the set $\mathbf I_n$ and the set of \emph{labelled Motzkin paths} of length $n$, namely, Motzkin paths whose down steps are labelled with an integer that does not exceed their height, while
the other steps are unlabelled. The set of labelled Motzkin paths of length $n$ will be denoted by $\mathcal L_n.$
Of course $\mathcal L_n\subset \mathcal H_n.$

In the present paper we need a bijection $\Psi$ that is a slightly modified version of the bijection $\Phi.$
The map $\Psi$ can be described as follows. Let $\pi \in \mathbf I_n.$ 

For every $i =
1, \ldots,n$:
\begin{itemize}
\item if $i$ is a fixed point for $\pi$ , draw a horizontal step;
\item if $i$ is the first element in a 2-cycle, draw an up step;
\item if $i$ is the second element in a 2-cycle $(j,i)$, draw a down step. Label this step with $h$, where $h$ is  the
number of cycles $(x,y)$ of $\pi$ such that $j<x<i<y$ .
\end{itemize}

For example, consider the involution $\pi=65382174$ whose standard cycle notation is $(7)(48)(3)(25)(16).$ Then 

\vspace{0.5cm}

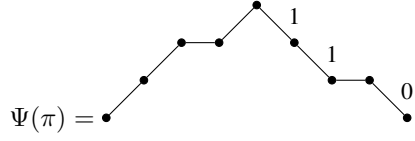
 \begin{displaymath}
\Psi(\pi)=\begin{tikzpicture}%[every node/.style={draw,shape=circle,minimum size=1mm, inner sep=0mm, fill=black}]
 \node[draw,shape=circle,minimum size=1mm, inner sep=0mm, fill=black] (8) at (8/2,0) {};
   \node[draw,shape=circle,minimum size=1mm, inner sep=0mm, fill=black] (7) at (7/2,1/2) {};
   \node[draw,shape=circle,minimum size=1mm, inner sep=0mm, fill=black] (6) at (6/2,1/2) {};
   \node[draw,shape=circle,minimum size=1mm, inner sep=0mm, fill=black] (5) at (5/2,1) {};
   \node[draw,shape=circle,minimum size=1mm, inner sep=0mm, fill=black] (4) at (4/2,3/2) {};
   \node[draw,shape=circle,minimum size=1mm, inner sep=0mm, fill=black] (3) at (3/2,1) {};
   \node[draw,shape=circle,minimum size=1mm, inner sep=0mm, fill=black] (2) at (2/2,1) {};
   \node[draw,shape=circle,minimum size=1mm, inner sep=0mm, fill=black] (1) at (1/2,1/2) {};
   \node[draw,shape=circle,minimum size=1mm, inner sep=0mm, fill=black] (0) at (0,0) {};
   \draw[] (0) -- (1) -- (2) -- (3) -- (4) -- (5) -- (6) -- (7) -- (8);
  \node[label={\small 1}] at (6/2, 1/2) {};
  \node[label={\small 1}] at (5/2, 2/2) {};
  \node[label={\small 0}] at (8/2, 0/2) {};
\end{tikzpicture}
\end{displaymath}
  
\vspace{0.5cm}

Our next aim is to prove the following result.
\begin{thm} 
\label{diagram}
The image under $\Gamma$ of the set $\mathbf S_n(\underline{132},\underline{123})$ is  $\mathcal L_n,$ and  the 
following diagram 
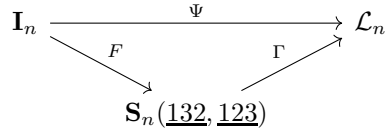
\begin{displaymath}\begin{tikzcd}
  \mathbf  I_n \arrow[rd,"F"] \arrow[rr, "\Psi"] & &  \mathcal L_n  \\
& \mathbf S_n(\underline{132},\underline{123}) \arrow[ru, "\Gamma"] &
  \end{tikzcd}\end{displaymath}
commutes.
\end{thm}
First of all, we characterize the image of the map $\Gamma$, when restricted to the set $\mathbf S_n(\underline{132},\underline{123})$.
\begin{prop}
Let $\pi\in\mathbf S_n$ and let $\Gamma(\pi)=(d,l).$
Then $\pi \in \mathbf S_n(\underline{132},\underline{123})$
if and only if
\begin{itemize}
\item the path $d=d_1\ldots d_n$ has no horizontal steps of color $c_2,$ and
\item for every index $i$, $l_i>0$ implies that $d_i$ is a down step.
\end{itemize}
\end{prop}
\proof

Firstly note that $\pi$ avoids the pattern $\underline{123}$ if and only if
the ascending runs of $\pi$ have length at most two. In this case the set of boarders of $\pi$ is empty and in $d$ there are no horizontal steps of color $c_2.$ 

Now let $\pi\in \mathbf S_n(\underline{132},\underline{123}).$ Suppose that there exists an integer $i$ such that $l_i>0$ and $d_i$ is either an up step or a horizontal step. 
By definition of the sequence $l,$ this implies that  the permutation $\pi$ contains three elements $\pi_s,\pi_{s+1},\pi_r$, with $\pi_r=i$, such that 
\begin{itemize}
\item $r>s+1$,
\item $\pi_s$ and $\pi_{s+1}$ are the head and the tail of an ascending run,
\item $\pi_s<\pi_r<\pi_{s+1}$, 
\item $\pi_r$ is either a head or a head-tail. 

\end{itemize} 
If $\pi_r$ is a head, then $\pi_{r+1}$ is the corresponding tail and $\pi_s,\pi_{r},\pi_{r+1}$ form an occurrence of
 $1\underline{23}.$
If $\pi_r$ is a head-tail, $\pi_{r-1}>\pi_r$ and $\pi_s,\pi_{r-1},\pi_r$ form an occurrence of $1\underline{32}.$
But this is impossible since, as noted above, $ \mathbf S_n(\underline{132},\underline{123})=\mathbf S_n(1\underline{32},1\underline{23}).$
On the other hand, if the permutation $\pi$ contains an occurrence of the pattern $\underline{132}$
corresponding to the elements $\pi_j,\pi_{j+1},\pi_{j+2},$  then $\pi_j<\pi_{j+2}<\pi_{j+1},$ $\pi_{j+2}$ is a head or a head-tail, and $\pi_j, \pi_{j+1}$ are the head and the tail of an ascending run. Hence $l_{\pi_{j+2}}>0.$
\endproof

Theorem \ref{diagram} now follows immediately from the description of the maps $F,$ $\Gamma$ and $\Psi$ and from the previous Proposition.

As an example consider the involution $\pi=(6)(48)(37)(2)(15)\in \mathbf I_8. $
Then the corresponding permutation in  $\mathbf  S_8(\underline{132},\underline{123})$ is  $F(\pi)=64837215$,  and
 \begin{displaymath}
\Psi(\pi)=
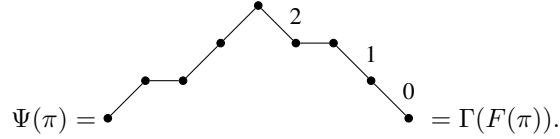
\begin{tikzpicture}%[every node/.style={draw,shape=circle,minimum size=1mm, inner sep=0mm, fill=black}]
 \node[draw,shape=circle,minimum size=1mm, inner sep=0mm, fill=black] (8) at (8/2,0) {};
   \node[draw,shape=circle,minimum size=1mm, inner sep=0mm, fill=black] (7) at (7/2,1/2) {};
   \node[draw,shape=circle,minimum size=1mm, inner sep=0mm, fill=black] (6) at (6/2,1) {};
   \node[draw,shape=circle,minimum size=1mm, inner sep=0mm, fill=black] (5) at (5/2,1) {};
   \node[draw,shape=circle,minimum size=1mm, inner sep=0mm, fill=black] (4) at (4/2,3/2) {};
   \node[draw,shape=circle,minimum size=1mm, inner sep=0mm, fill=black] (3) at (3/2,1) {};
   \node[draw,shape=circle,minimum size=1mm, inner sep=0mm, fill=black] (2) at (2/2,1/2) {};
   \node[draw,shape=circle,minimum size=1mm, inner sep=0mm, fill=black] (1) at (1/2,1/2) {};
   \node[draw,shape=circle,minimum size=1mm, inner sep=0mm, fill=black] (0) at (0,0) {};
   \draw[] (0) -- (1) -- (2) -- (3) -- (4) -- (5) -- (6) -- (7) -- (8);
  \node[label={\small 1}] at (7/2, 1/2) {};
  \node[label={\small 2}] at (5/2, 2/2) {};
  \node[label={\small 0}] at (8/2, 0/2) {};
\end{tikzpicture}=\Gamma(F(\pi)).
\end{displaymath}

The subset of $\mathcal L_n$ of labelled Motzkin paths of length $n$ all of whose labels are zero is obviously isomorphic to the set $\mathcal M_n$ of Motzkin paths. It is possible to characterize the preimage of this set under the maps $\Gamma$ and $\Psi$ in terms of pattern avoiding permutations. In fact, in \cite{BBCS} the following result is proved.\begin{prop}\label{132zero}
Let $\pi\in \mathbf S_n$ and let $\Gamma(\pi)=(d,l).$ Then
\begin{equation} \pi\in \mathbf S_n(132)\mbox{    if and only if   } l=(0,\ldots ,0).\end{equation}
\end{prop}
As a consequence $\Gamma$ induces a bijection between $S_n({132},\underline{123})$ and $\mathcal M_n.$

Moreover, in \cite[Theorem 9]{Ba6} it has been shown that the map  $\Psi$ induces a bijection between the set of involutions avoiding the pattern $3412$ and the set of labelled Motzkin paths of length $n$ all of whose labels are zero.
These results imply that the following diagram  
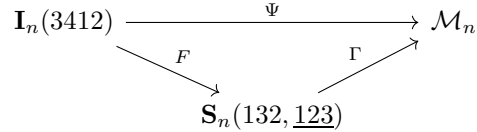
\begin{displaymath}
\begin{tikzcd}
  \mathbf I_n(3412)  \arrow[rd,"F"] \arrow[rr, "\Psi"] & &  \mathcal M_n  \\
& \mathbf S_n({132},\underline{123}) \arrow[ru, "\Gamma"] &
\end{tikzcd}
\end{displaymath}
commutes. 

In the following sections we show how some statistics over the sets $\mathbf  S_n({132},\underline{123})$ and $\mathbf I_n(3412)$ can be translated into statistics over Motzkin paths.

\section{Inversions and descents over $\mathbf I_n(3412)$}

Let $d$ be a Motzkin path. A \emph{tunnel} in $d$ is a horizontal segment between two lattice points of $d$ lying weakly below $d$ and containing exactly two lattice points of $d.$ 
Note that each horizontal step of $d$ is a tunnel. We will call the horizontal steps \emph{trivial tunnels.}

We recall that each non-empty Motzkin  path $m$ can be decomposed either as $Hm',$ where $m'$ is an arbitrary
Motzkin path, or as $Um'Dm'',$ where $m'$ and $m''$ are arbitrary Motzkin paths. 
This decomposition is called \emph{first return decomposition}.
The definition of the map $\Psi$ implies that each $2$-cycle of an involution $\pi\in \mathbf I_n(3412)$  corresponds to a  non-trivial tunnel of $\Psi(\pi)$ and vice-versa.

Let $\pi=\pi_1\ldots \pi_n \in \mathbf S_n.$ An \emph{inversion} of $\pi$ is a pair $(i,j)$ with $i<j$ such that $\pi_i>\pi_j.$
In this case we will say that the symbol $\pi_i$ \emph{is in inversion} with the symbol $\pi_j.$ 
The number of inversions of the permutation $\pi$ will be denoted by $\textrm{inv}(\pi).$

The permutation $\pi$ has a \emph{descent} at position $i$ if $\pi_i>\pi_{i+1}.$ Otherwise, $\pi$ has an \emph{ascent} at position $i$.The number of descents of $\pi$ will be denoted by $\textrm{des}(\pi).$ 

Now we want to study the joint distribution of the statistics inv, des, fix over the set 
\begin{displaymath} \mathbf I(3412):=\bigcup_{n\geq 0} \mathbf I_n(3412),\end{displaymath} where $\textrm{fix}(\pi)$ denotes the number of fixed points of $\pi,$
namely, determine an expression for the generating function
\begin{equation} \label{funzione}
F(x,y,z,w)=\sum_{n\geq 0}\sum_{\pi\in\mathbf I_n(3412)} x^n y^{\textrm{inv}(\pi)}z^{\textrm{des}(\pi)}w^{\textrm{fix}(\pi)}.
 \end{equation}

First of all we prove a preliminary result. 
\begin{lem}
Let $\pi\in\mathbf I_n(3412).$  Then
\begin{equation}\textrm{inv}(\pi)=2A-t,\end{equation}
where $A$ is the area between the path $\Psi(\pi)$ and the $x$-axis and $t$ is the number of non-trivial tunnels of $\Psi(\pi).$
\end{lem}
\proof
Write $\pi$ as $\pi_1\pi_2\ldots \pi_n.$ Let $i$ be the least index such that $\pi_i>i.$ Then $(i,\pi_i)$ is a cycle of $\pi.$ Hence the symbol $\pi_i$ is in inversion with all the symbols
$\pi_{i+1},\pi_{i+2}\ldots \pi_{i+k}=i,$ where $k=\pi_i-i.$
In fact, suppose by contradiction that there exist an index $r,$ with $1\leq r\leq k-1,$ such that  $\pi_{i+r}>\pi_i,$ then $\pi_i,\pi_{i+r},i,r.$ would be an occurrence of the pattern $3412.$

For the same reason the symbol $i$ is in inversion with the $k-1$ elements $\pi_{i+1},\ldots,\pi_{i+k-1}.$ Here we excluded the inversion $(i,\pi_i).$

On the other hand, let $T$ be the tunnel in the Motzkin path $\Psi(\pi)$ corresponding to the cycle $(i,\pi_i).$ The area of the trapezoid with height one and  $T$ as a basis
is precisely $k=\pi_i-i.$

Repeating the preceding argument on the involution obtained from $\pi$ by deleting the symbols $i$ and $\pi_i$ we get the assertion.
\endproof

\begin{lem}
\label{lem42}
Let $\pi\in\mathbf I_n(3412).$ 
The descents of $\pi$ correspond bijectively to the occurrences in $\Psi(\pi)$ of the following subwords:
$UU,$ $DD,$ $UH,$ $HD$ and $UD.$ 
\end{lem}
\proof
Suppose that one of these subwords occurs in the Motzkin path $\Psi(\pi).$ Let $v$ be this subword and $i$ be the position of the first step of $v.$
If $v=UU,$ then both $\pi_i$ and $\pi_{i+1}$ are the greater elements in their respective $2$-cycles and hence $\pi_i>\pi_{i+1}.$
If $v=UH,$ then $\pi_i$ is the greater element in its $2$-cycle while $\pi_{i+1}$ is a fixed point and hence $\pi_i>\pi_{i+1}.$
The other cases can be treated in a similar way. 
\endproof

We recall that a \emph{weak valley} of a Motzkin path is an occurrence of one of  the following consecutive patterns:
\begin{displaymath} HH, \quad HU, \quad DH, \quad DU.\end{displaymath}
The preceding result yields immediately:
\begin{cor}
The distribution of ascents over $n$-involutions avoiding the pattern 3412 is the same as the distribution of weak valleys over Motzkin paths of length $n$.
\end{cor}
This implies that the generating function $G(x,z)$ of Motzkin paths according to the length ($x$) and the number of weak valleys  ($z$) can be deduced from the function $F(x,y,z,w)$ appearing in Formula (\ref{funzione}) as follows:
\begin{equation}
G(x,z)=1+\frac{F(xz,1,1/z,1)-1}{z}.
\end{equation}
Similarly, $1+\frac{F(xz,1,1/z,0)-1}{z}$ gives the generating function of Dyck paths according to the length and the number of valleys, which is essentially given by Narayana polynomials.

The above lemmas imply that the generating function $F$ satisfies the following equation obtained by the first return decomposition for Motzkin paths.

\begin{equation}\label{ricorrenzagenfun}
\begin{split}
F(x,y,z,w) & =1+xwF(x,y,z,w)+x^2yzF(x,y,z,w) \\
& + x^2yz^2(F(xy^2,y,z,w)-1)\cdot F(x,y,z,w). 
\end{split}
\end{equation}

In fact, the terms in the right hand side of the previous equation correspond to Motzkin paths either  empty or  of the form $Hm,$ $UDm,$ $Um'Dm,$
with $m'$ non-empty, respectively. 

From equation (\ref{ricorrenzagenfun}) we get easily the following continued fraction expression for $F$.

\begin{thm}
\begin{equation*}
\begin{split}
F(x,y,z,w)= &  \frac{1}{1-xw-x^2yz+x^2yz^2-x^2yz^2F(xy^2,y,z,w)}=\\
 & =\cfrac{1}{1+b_0-\cfrac{c_0}{1+b_1-\cfrac{c_1}{1+b_2-\cfrac{c_2}{1+\ldots}}}}
\end{split}
\end{equation*}
where $b_i=-xy^{2i}w-x^2y^{2i+1}z+x^2y^{2i+1}z^2 $
and $c_i=x^2y^{2i+1}z^2,$ $i\geq 0.$
 \end{thm}

\section{The distribution of consecutive patterns in $\mathbf I_n(3412)$}\label{involutionspatterns}

Lemma \ref{lem42} allows us to translate every three-letter subword of a Motzkin path into an occurrence of a consecutive pattern of the corresponding involution in $\mathbf I_n(3412)$.

\begin{thm}\label{corrispondenze}
 Let $\pi\in\mathbf I_n(3412)$ and let $\Psi(\pi)$ be the corresponding
 Motzkin path. Then a subword of $\Psi(\pi)$ of length three corresponds to an occurrence of a  consecutive pattern in $\pi$
 according to the following table. 
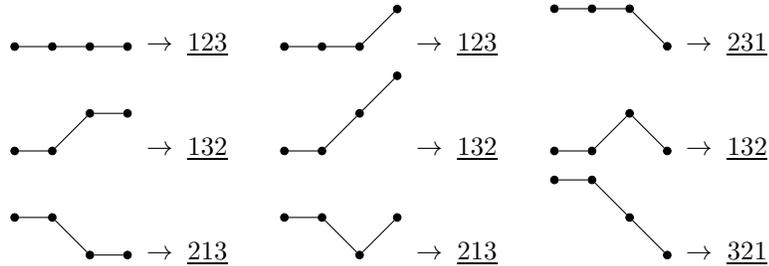
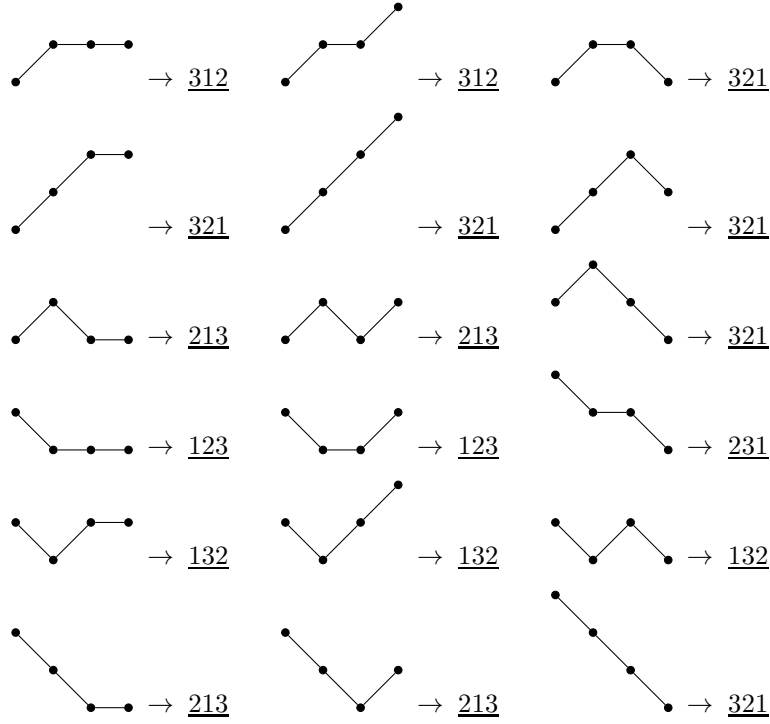
\begin{displaymath}
\begin{tikzpicture}%[every node/.style={draw,shape=circle,minimum size=1mm, inner sep=0mm, fill=black}]
   \node[draw,shape=circle,minimum size=1mm, inner sep=0mm, fill=black] (3) at (3/2,0) {};
   \node[draw,shape=circle,minimum size=1mm, inner sep=0mm, fill=black] (2) at (2/2,0) {};
   \node[draw,shape=circle,minimum size=1mm, inner sep=0mm, fill=black] (1) at (1/2,0) {};
   \node[draw,shape=circle,minimum size=1mm, inner sep=0mm, fill=black] (0) at (0,0) {};
   \draw[] (0) -- (1) -- (2) -- (3);
\end{tikzpicture}\; \rightarrow  \; \underline{123}\qquad
\begin{tikzpicture}%[every node/.style={draw,shape=circle,minimum size=1mm, inner sep=0mm, fill=black}]
   \node[draw,shape=circle,minimum size=1mm, inner sep=0mm, fill=black] (3) at (3/2,1/2) {};
   \node[draw,shape=circle,minimum size=1mm, inner sep=0mm, fill=black] (2) at (2/2,0) {};
   \node[draw,shape=circle,minimum size=1mm, inner sep=0mm, fill=black] (1) at (1/2,0) {};
   \node[draw,shape=circle,minimum size=1mm, inner sep=0mm, fill=black] (0) at (0,0) {};
   \draw[] (0) -- (1) -- (2) -- (3);
\end{tikzpicture}\; \rightarrow  \; \underline{123}\qquad
\begin{tikzpicture}%[every node/.style={draw,shape=circle,minimum size=1mm, inner sep=0mm, fill=black}]
   \node[draw,shape=circle,minimum size=1mm, inner sep=0mm, fill=black] (3) at (3/2,-1/2) {};
   \node[draw,shape=circle,minimum size=1mm, inner sep=0mm, fill=black] (2) at (2/2,0) {};
   \node[draw,shape=circle,minimum size=1mm, inner sep=0mm, fill=black] (1) at (1/2,0) {};
   \node[draw,shape=circle,minimum size=1mm, inner sep=0mm, fill=black] (0) at (0,0) {};
   \draw[] (0) -- (1) -- (2) -- (3);
\end{tikzpicture}\; \rightarrow  \; \underline{231}\end{displaymath}
\begin{displaymath}
\begin{tikzpicture}%[every node/.style={draw,shape=circle,minimum size=1mm, inner sep=0mm, fill=black}]
   \node[draw,shape=circle,minimum size=1mm, inner sep=0mm, fill=black] (3) at (3/2,1/2) {};
   \node[draw,shape=circle,minimum size=1mm, inner sep=0mm, fill=black] (2) at (2/2,1/2) {};
   \node[draw,shape=circle,minimum size=1mm, inner sep=0mm, fill=black] (1) at (1/2,0) {};
   \node[draw,shape=circle,minimum size=1mm, inner sep=0mm, fill=black] (0) at (0,0) {};
   \draw[] (0) -- (1) -- (2) -- (3);
\end{tikzpicture}\; \rightarrow  \; \underline{132}\qquad
\begin{tikzpicture}%[every node/.style={draw,shape=circle,minimum size=1mm, inner sep=0mm, fill=black}]
   \node[draw,shape=circle,minimum size=1mm, inner sep=0mm, fill=black] (3) at (3/2,1) {};
   \node[draw,shape=circle,minimum size=1mm, inner sep=0mm, fill=black] (2) at (2/2,1/2) {};
   \node[draw,shape=circle,minimum size=1mm, inner sep=0mm, fill=black] (1) at (1/2,0) {};
   \node[draw,shape=circle,minimum size=1mm, inner sep=0mm, fill=black] (0) at (0,0) {};
   \draw[] (0) -- (1) -- (2) -- (3);
\end{tikzpicture}\; \rightarrow  \; \underline{132}\qquad
\begin{tikzpicture}%[every node/.style={draw,shape=circle,minimum size=1mm, inner sep=0mm, fill=black}]
   \node[draw,shape=circle,minimum size=1mm, inner sep=0mm, fill=black] (3) at (3/2,0) {};
   \node[draw,shape=circle,minimum size=1mm, inner sep=0mm, fill=black] (2) at (2/2,1/2) {};
   \node[draw,shape=circle,minimum size=1mm, inner sep=0mm, fill=black] (1) at (1/2,0) {};
   \node[draw,shape=circle,minimum size=1mm, inner sep=0mm, fill=black] (0) at (0,0) {};
   \draw[] (0) -- (1) -- (2) -- (3);
\end{tikzpicture}\; \rightarrow  \; \underline{132}\end{displaymath}
\begin{displaymath}
\begin{tikzpicture}%[every node/.style={draw,shape=circle,minimum size=1mm, inner sep=0mm, fill=black}]
   \node[draw,shape=circle,minimum size=1mm, inner sep=0mm, fill=black] (3) at (3/2,-1/2) {};
   \node[draw,shape=circle,minimum size=1mm, inner sep=0mm, fill=black] (2) at (2/2,-1/2) {};
   \node[draw,shape=circle,minimum size=1mm, inner sep=0mm, fill=black] (1) at (1/2,0) {};
   \node[draw,shape=circle,minimum size=1mm, inner sep=0mm, fill=black] (0) at (0,0) {};
   \draw[] (0) -- (1) -- (2) -- (3);
\end{tikzpicture}\; \rightarrow  \; \underline{213}\qquad
\begin{tikzpicture}%[every node/.style={draw,shape=circle,minimum size=1mm, inner sep=0mm, fill=black}]
   \node[draw,shape=circle,minimum size=1mm, inner sep=0mm, fill=black] (3) at (3/2,0) {};
   \node[draw,shape=circle,minimum size=1mm, inner sep=0mm, fill=black] (2) at (2/2,-1/2) {};
   \node[draw,shape=circle,minimum size=1mm, inner sep=0mm, fill=black] (1) at (1/2,0) {};
   \node[draw,shape=circle,minimum size=1mm, inner sep=0mm, fill=black] (0) at (0,0) {};
   \draw[] (0) -- (1) -- (2) -- (3);
\end{tikzpicture}\; \rightarrow  \; \underline{213}\qquad
\begin{tikzpicture}%[every node/.style={draw,shape=circle,minimum size=1mm, inner sep=0mm, fill=black}]
   \node[draw,shape=circle,minimum size=1mm, inner sep=0mm, fill=black] (3) at (3/2,-1) {};
   \node[draw,shape=circle,minimum size=1mm, inner sep=0mm, fill=black] (2) at (2/2,-1/2) {};
   \node[draw,shape=circle,minimum size=1mm, inner sep=0mm, fill=black] (1) at (1/2,0) {};
   \node[draw,shape=circle,minimum size=1mm, inner sep=0mm, fill=black] (0) at (0,0) {};
   \draw[] (0) -- (1) -- (2) -- (3);
\end{tikzpicture}\; \rightarrow  \; \underline{321}\end{displaymath}
\begin{displaymath}
\begin{tikzpicture}%[every node/.style={draw,shape=circle,minimum size=1mm, inner sep=0mm, fill=black}]
   \node[draw,shape=circle,minimum size=1mm, inner sep=0mm, fill=black] (3) at (3/2,1/2) {};
   \node[draw,shape=circle,minimum size=1mm, inner sep=0mm, fill=black] (2) at (2/2,1/2) {};
   \node[draw,shape=circle,minimum size=1mm, inner sep=0mm, fill=black] (1) at (1/2,1/2) {};
   \node[draw,shape=circle,minimum size=1mm, inner sep=0mm, fill=black] (0) at (0,0) {};
   \draw[] (0) -- (1) -- (2) -- (3);
\end{tikzpicture}\; \rightarrow  \; \underline{312}\qquad
\begin{tikzpicture}%[every node/.style={draw,shape=circle,minimum size=1mm, inner sep=0mm, fill=black}]
   \node[draw,shape=circle,minimum size=1mm, inner sep=0mm, fill=black] (3) at (3/2,1) {};
   \node[draw,shape=circle,minimum size=1mm, inner sep=0mm, fill=black] (2) at (2/2,1/2) {};
   \node[draw,shape=circle,minimum size=1mm, inner sep=0mm, fill=black] (1) at (1/2,1/2) {};
   \node[draw,shape=circle,minimum size=1mm, inner sep=0mm, fill=black] (0) at (0,0) {};
   \draw[] (0) -- (1) -- (2) -- (3);
\end{tikzpicture}\; \rightarrow  \; \underline{312}\qquad
\begin{tikzpicture}%[every node/.style={draw,shape=circle,minimum size=1mm, inner sep=0mm, fill=black}]
   \node[draw,shape=circle,minimum size=1mm, inner sep=0mm, fill=black] (3) at (3/2,0) {};
   \node[draw,shape=circle,minimum size=1mm, inner sep=0mm, fill=black] (2) at (2/2,1/2) {};
   \node[draw,shape=circle,minimum size=1mm, inner sep=0mm, fill=black] (1) at (1/2,1/2) {};
   \node[draw,shape=circle,minimum size=1mm, inner sep=0mm, fill=black] (0) at (0,0) {};
   \draw[] (0) -- (1) -- (2) -- (3);
\end{tikzpicture}\; \rightarrow  \; \underline{321}\end{displaymath}
\begin{displaymath}
\begin{tikzpicture}%[every node/.style={draw,shape=circle,minimum size=1mm, inner sep=0mm, fill=black}]
   \node[draw,shape=circle,minimum size=1mm, inner sep=0mm, fill=black] (3) at (3/2,1) {};
   \node[draw,shape=circle,minimum size=1mm, inner sep=0mm, fill=black] (2) at (2/2,1) {};
   \node[draw,shape=circle,minimum size=1mm, inner sep=0mm, fill=black] (1) at (1/2,1/2) {};
   \node[draw,shape=circle,minimum size=1mm, inner sep=0mm, fill=black] (0) at (0,0) {};
   \draw[] (0) -- (1) -- (2) -- (3);
\end{tikzpicture}\; \rightarrow  \; \underline{321}\qquad
\begin{tikzpicture}%[every node/.style={draw,shape=circle,minimum size=1mm, inner sep=0mm, fill=black}]
   \node[draw,shape=circle,minimum size=1mm, inner sep=0mm, fill=black] (3) at (3/2,3/2) {};
   \node[draw,shape=circle,minimum size=1mm, inner sep=0mm, fill=black] (2) at (2/2,1) {};
   \node[draw,shape=circle,minimum size=1mm, inner sep=0mm, fill=black] (1) at (1/2,1/2) {};
   \node[draw,shape=circle,minimum size=1mm, inner sep=0mm, fill=black] (0) at (0,0) {};
   \draw[] (0) -- (1) -- (2) -- (3);
\end{tikzpicture}\; \rightarrow  \; \underline{321}\qquad
\begin{tikzpicture}%[every node/.style={draw,shape=circle,minimum size=1mm, inner sep=0mm, fill=black}]
   \node[draw,shape=circle,minimum size=1mm, inner sep=0mm, fill=black] (3) at (3/2,1/2) {};
   \node[draw,shape=circle,minimum size=1mm, inner sep=0mm, fill=black] (2) at (2/2,1) {};
   \node[draw,shape=circle,minimum size=1mm, inner sep=0mm, fill=black] (1) at (1/2,1/2) {};
   \node[draw,shape=circle,minimum size=1mm, inner sep=0mm, fill=black] (0) at (0,0) {};
   \draw[] (0) -- (1) -- (2) -- (3);
\end{tikzpicture}\; \rightarrow  \; \underline{321}\end{displaymath}
\begin{displaymath}
\begin{tikzpicture}%[every node/.style={draw,shape=circle,minimum size=1mm, inner sep=0mm, fill=black}]
   \node[draw,shape=circle,minimum size=1mm, inner sep=0mm, fill=black] (3) at (3/2,0) {};
   \node[draw,shape=circle,minimum size=1mm, inner sep=0mm, fill=black] (2) at (2/2,0) {};
   \node[draw,shape=circle,minimum size=1mm, inner sep=0mm, fill=black] (1) at (1/2,1/2) {};
   \node[draw,shape=circle,minimum size=1mm, inner sep=0mm, fill=black] (0) at (0,0) {};
   \draw[] (0) -- (1) -- (2) -- (3);
\end{tikzpicture}\; \rightarrow  \; \underline{213}\qquad
\begin{tikzpicture}%[every node/.style={draw,shape=circle,minimum size=1mm, inner sep=0mm, fill=black}]
   \node[draw,shape=circle,minimum size=1mm, inner sep=0mm, fill=black] (3) at (3/2,1/2) {};
   \node[draw,shape=circle,minimum size=1mm, inner sep=0mm, fill=black] (2) at (2/2,0) {};
   \node[draw,shape=circle,minimum size=1mm, inner sep=0mm, fill=black] (1) at (1/2,1/2) {};
   \node[draw,shape=circle,minimum size=1mm, inner sep=0mm, fill=black] (0) at (0,0) {};
   \draw[] (0) -- (1) -- (2) -- (3);
\end{tikzpicture}\; \rightarrow  \; \underline{213}\qquad
\begin{tikzpicture}%[every node/.style={draw,shape=circle,minimum size=1mm, inner sep=0mm, fill=black}]
   \node[draw,shape=circle,minimum size=1mm, inner sep=0mm, fill=black] (3) at (3/2,-1/2) {};
   \node[draw,shape=circle,minimum size=1mm, inner sep=0mm, fill=black] (2) at (2/2,0) {};
   \node[draw,shape=circle,minimum size=1mm, inner sep=0mm, fill=black] (1) at (1/2,1/2) {};
   \node[draw,shape=circle,minimum size=1mm, inner sep=0mm, fill=black] (0) at (0,0) {};
   \draw[] (0) -- (1) -- (2) -- (3);
\end{tikzpicture}\; \rightarrow  \; \underline{321}\end{displaymath}
\begin{displaymath}
\begin{tikzpicture}%[every node/.style={draw,shape=circle,minimum size=1mm, inner sep=0mm, fill=black}]
   \node[draw,shape=circle,minimum size=1mm, inner sep=0mm, fill=black] (3) at (3/2,-1/2) {};
   \node[draw,shape=circle,minimum size=1mm, inner sep=0mm, fill=black] (2) at (2/2,-1/2) {};
   \node[draw,shape=circle,minimum size=1mm, inner sep=0mm, fill=black] (1) at (1/2,-1/2) {};
   \node[draw,shape=circle,minimum size=1mm, inner sep=0mm, fill=black] (0) at (0,0) {};
   \draw[] (0) -- (1) -- (2) -- (3);
\end{tikzpicture}\; \rightarrow  \; \underline{123}\qquad
\begin{tikzpicture}%[every node/.style={draw,shape=circle,minimum size=1mm, inner sep=0mm, fill=black}]
   \node[draw,shape=circle,minimum size=1mm, inner sep=0mm, fill=black] (3) at (3/2,0) {};
   \node[draw,shape=circle,minimum size=1mm, inner sep=0mm, fill=black] (2) at (2/2,-1/2) {};
   \node[draw,shape=circle,minimum size=1mm, inner sep=0mm, fill=black] (1) at (1/2,-1/2) {};
   \node[draw,shape=circle,minimum size=1mm, inner sep=0mm, fill=black] (0) at (0,0) {};
   \draw[] (0) -- (1) -- (2) -- (3);
\end{tikzpicture}\; \rightarrow  \; \underline{123}\qquad
\begin{tikzpicture}%[every node/.style={draw,shape=circle,minimum size=1mm, inner sep=0mm, fill=black}]
   \node[draw,shape=circle,minimum size=1mm, inner sep=0mm, fill=black] (3) at (3/2,-1) {};
   \node[draw,shape=circle,minimum size=1mm, inner sep=0mm, fill=black] (2) at (2/2,-1/2) {};
   \node[draw,shape=circle,minimum size=1mm, inner sep=0mm, fill=black] (1) at (1/2,-1/2) {};
   \node[draw,shape=circle,minimum size=1mm, inner sep=0mm, fill=black] (0) at (0,0) {};
   \draw[] (0) -- (1) -- (2) -- (3);
\end{tikzpicture}\; \rightarrow  \; \underline{231}\end{displaymath}
\begin{displaymath}
\begin{tikzpicture}%[every node/.style={draw,shape=circle,minimum size=1mm, inner sep=0mm, fill=black}]
   \node[draw,shape=circle,minimum size=1mm, inner sep=0mm, fill=black] (3) at (3/2,0) {};
   \node[draw,shape=circle,minimum size=1mm, inner sep=0mm, fill=black] (2) at (2/2,0) {};
   \node[draw,shape=circle,minimum size=1mm, inner sep=0mm, fill=black] (1) at (1/2,-1/2) {};
   \node[draw,shape=circle,minimum size=1mm, inner sep=0mm, fill=black] (0) at (0,0) {};
   \draw[] (0) -- (1) -- (2) -- (3);
\end{tikzpicture}\; \rightarrow  \; \underline{132}\qquad
\begin{tikzpicture}%[every node/.style={draw,shape=circle,minimum size=1mm, inner sep=0mm, fill=black}]
   \node[draw,shape=circle,minimum size=1mm, inner sep=0mm, fill=black] (3) at (3/2,1/2) {};
   \node[draw,shape=circle,minimum size=1mm, inner sep=0mm, fill=black] (2) at (2/2,0) {};
   \node[draw,shape=circle,minimum size=1mm, inner sep=0mm, fill=black] (1) at (1/2,-1/2) {};
   \node[draw,shape=circle,minimum size=1mm, inner sep=0mm, fill=black] (0) at (0,0) {};
   \draw[] (0) -- (1) -- (2) -- (3);
\end{tikzpicture}\; \rightarrow  \; \underline{132}\qquad
\begin{tikzpicture}%[every node/.style={draw,shape=circle,minimum size=1mm, inner sep=0mm, fill=black}]
   \node[draw,shape=circle,minimum size=1mm, inner sep=0mm, fill=black] (3) at (3/2,-1/2) {};
   \node[draw,shape=circle,minimum size=1mm, inner sep=0mm, fill=black] (2) at (2/2,0) {};
   \node[draw,shape=circle,minimum size=1mm, inner sep=0mm, fill=black] (1) at (1/2,-1/2) {};
   \node[draw,shape=circle,minimum size=1mm, inner sep=0mm, fill=black] (0) at (0,0) {};
   \draw[] (0) -- (1) -- (2) -- (3);
\end{tikzpicture}\; \rightarrow  \; \underline{132}\end{displaymath}
\begin{displaymath}
\begin{tikzpicture}%[every node/.style={draw,shape=circle,minimum size=1mm, inner sep=0mm, fill=black}]
   \node[draw,shape=circle,minimum size=1mm, inner sep=0mm, fill=black] (3) at (3/2,-1) {};
   \node[draw,shape=circle,minimum size=1mm, inner sep=0mm, fill=black] (2) at (2/2,-1) {};
   \node[draw,shape=circle,minimum size=1mm, inner sep=0mm, fill=black] (1) at (1/2,-1/2) {};
   \node[draw,shape=circle,minimum size=1mm, inner sep=0mm, fill=black] (0) at (0,0) {};
   \draw[] (0) -- (1) -- (2) -- (3);
\end{tikzpicture}\; \rightarrow  \; \underline{213}\qquad
\begin{tikzpicture}%[every node/.style={draw,shape=circle,minimum size=1mm, inner sep=0mm, fill=black}]
   \node[draw,shape=circle,minimum size=1mm, inner sep=0mm, fill=black] (3) at (3/2,-1/2) {};
   \node[draw,shape=circle,minimum size=1mm, inner sep=0mm, fill=black] (2) at (2/2,-1) {};
   \node[draw,shape=circle,minimum size=1mm, inner sep=0mm, fill=black] (1) at (1/2,-1/2) {};
   \node[draw,shape=circle,minimum size=1mm, inner sep=0mm, fill=black] (0) at (0,0) {};
   \draw[] (0) -- (1) -- (2) -- (3);
\end{tikzpicture}\; \rightarrow  \; \underline{213}\qquad
\begin{tikzpicture}%[every node/.style={draw,shape=circle,minimum size=1mm, inner sep=0mm, fill=black}]
   \node[draw,shape=circle,minimum size=1mm, inner sep=0mm, fill=black] (3) at (3/2,-3/2) {};
   \node[draw,shape=circle,minimum size=1mm, inner sep=0mm, fill=black] (2) at (2/2,-1) {};
   \node[draw,shape=circle,minimum size=1mm, inner sep=0mm, fill=black] (1) at (1/2,-1/2) {};
   \node[draw,shape=circle,minimum size=1mm, inner sep=0mm, fill=black] (0) at (0,0) {};
   \draw[] (0) -- (1) -- (2) -- (3);
\end{tikzpicture}\; \rightarrow  \; \underline{321}\\
\end{displaymath}
\end{thm}

Now we enumerate the involutions in $\mathbf I_n(3412)$ according to the occurrences of a given consecutive pattern of length three and the
number of fixed points. 

First of all, observe that an involution $\pi\in \mathbf I_n(3412)$ has $k$ occurrences of $\underline{213}$ ($\underline{312}$) and $f$ fixed points if and only if $RC(\pi)$
has $k$ occurrences of $\underline{132}$ ($\underline{231}$,  respectively) and $f$ fixed points, where $RC(\pi)$ is the \emph{reverse-complement} of the permutation $\pi,$
namely,  $RC(\pi=\pi_1\ldots \pi_n)=n+1-\pi_n\ldots n+1-\pi_1.$

Hence we can restrict our attention to the consecutive patterns $\underline{123}, $ $\underline{321},$
$\underline{132}$ and $\underline{231}$.

\subsection{The pattern \underline{123}}

Let $a_{n,k,f}$ be the number of involutions $\pi\in\mathbf I_n(3412)$ with $k$ occurrences of $\underline{123}$ and with $f$
fixed points. 

Our goal is to find a formula for the generating function \begin{equation} F=\sum a_{n,k,f}x^nt^kz^f.\end{equation}
To this aim we use a variation of the Goulden-Jackson cluster method (see \cite[p. 128]{Goulden}).
In the Appendix we give a detailed description of the notations and the results that will be used.

By Theorem \ref{corrispondenze}, the pattern $\underline{123}$ in $\pi$ corresponds to occurrences in $\Psi(\pi)$ of subwords in the set
  $S=\{H^3, H^2U, DH^2, DHU\}.$
These subwords give rise to the clusters of type $H^j$ with $j\geq 3,$ $H^jU$ with $j\geq 2,$
$DH^j$ with $j\geq 2,$ and  $DH^jU $ with $j\geq 1.$

Note that
\begin{itemize}
\item the cluster $H^j$ reduces to a horizontal step and has depth $0,$ 
\item the cluster
$H^jU$ reduces to an up step and has depth $0,$ 
\item the cluster  $DH^j$ reduces to a down step and has depth $0,$ 
\item the cluster $DH^jU$ reduces to a horizontal step and has depth $-1.$ 
\end{itemize}
To find $F$ we will use Theorem \ref{teoremacluster} of the Appendix.

First of all, we determine the generating functions $A_H(x,t,z),$ $A'_H(x,t,z),$ $A_D(x,t,z),$   $A'_D(x,t,z),$  $A_U(x,t,z)$ and  $A'_U(x,t,z).$

Consider the cluster $H^j, $ $j\geq 5.$  This can be obtained by either  juxtaposing a horizontal step to the right of $H^{j-1}$
and adding an occurrence of the subword $H^3$  that covers the last two letters of $H^{j-1}$
\begin{equation} H^j=\lefteqn{\underbrace{\phantom{HH\ldots HH}}_{H^{j-1}}}HH\ldots\overbrace{HHH}^{H^3}\end{equation}
or  juxtaposing two horizontal steps to the right of $H^{j-2}$
and adding an occurrence of the subword $H^3$ that covers the last letter of $H^{j-2}$
\begin{equation}H^j=\lefteqn{\underbrace{\phantom{HH\ldots H}}_{H^{j-2}}}HH\ldots\overbrace{HHH}^{H^3}.\end{equation}
Note that in these two cases the number of occurrences of $H^3$ increases by one, the number of horizontal steps and the length increase by one in the first case and by two in the second case.
As a consequence, the generating function for clusters of this kind is
\begin{equation}\frac{x^3z^3t}{1-xzt-x^2z^2t},\end{equation}
where the variables $x,t,$ and $z$ represent length, number of occurrences of the subwords in $S$ and number of $H,$ respectively.

Similarly, the cluster $DH^j,$ $j\geq 4,$  can be obtained
in the two ways depicted below
\begin{equation} DH^j=\lefteqn{\underbrace{\phantom{DH\ldots HH}}_{DH^{j-1}}}DH\ldots\overbrace{HHH}^{H^3}\end{equation}
or
\begin{equation} DH^j=\lefteqn{\underbrace{\phantom{DH\ldots H}}_{DH^{j-2}}}DH\ldots\overbrace{HHH}^{H^3}\end{equation}
hence, the corresponding generating function is
\begin{equation} \label{DHH}
\frac{x^3z^2t}{1-xzt-x^2z^2t}.
\end{equation}

By similar arguments the generating function for the cluster $H^jU$ is
\begin{equation} \frac{x^3z^2t}{1-xzt-x^2z^2t}.\end{equation}

Lastly, the cluster $DH^jU,$ $j\geq 3$ can be obtained by either  juxtaposing  the  letter $U$ to the right of $DH^{j}$ and adding
an occurrence of $H^2U$ that covers the last two  letters of $DH^{j-1},$ or juxtaposing the  letters $HU$ to the right of $DH^{j-1}$
and adding an occurrence of $H^2U$ that covers the last letter of $DH^{j-1}.$ 
By formula (\ref{DHH}) we get the following expression for the generating function of the cluster of the form $DH^jU,$ $j\geq 2$  
\begin{equation} \frac{x^3z^2t}{1-xzt-x^2z^2t}\cdot (xt+x^2tz)\end{equation}

The cluster $DHU$ must be considered separately.
Its contribution is $x^3tz.$

As a consequence we have
\begin{equation} A_H(x,t,k)=\frac{x^3z^2t(z+xt+x^2tz)}{1-xzt-x^2z^2t}+x^3tz,\end{equation}
\begin{equation} A'_H(x,t,z)=\frac{x^3z^3t}{1-xzt-x^2z^2t}\end{equation} and
\begin{equation} A_D(x,t,z)=A'_D(x,t,z)=A_U(x,t,z)=A'_U(x,t,z)=\frac{x^3z^2t}{1-xzt-x^2z^2t}.\end{equation}

Now we are in position to apply Theorem \ref{teoremacluster}, hence
finding
the generating function $F_{\underline{123}}$ evaluated in $x,t+1,z.$ After the substitution 
 $t\leftarrow t-1$  we get the following expression for $F_{\underline{123}}(x,t,z).$
 \begin{thm}
\begin{equation}F_{\underline{123}}(x,t,z)=\frac{2A^2}{2(1-B)A^2-A^2+A^2C+A^2\sqrt{(1-C)^2-4A^2}}\end{equation}
where
\begin{equation}A=x+\frac{x^3z^2(t-1)}{1-xz(t-1)-x^2z^2(t-1)}, \end{equation}
\begin{equation}B=xz+\frac{x^3z^3(t-1)}{1-xz(t-1)-x^2z^2(t-1)}\end{equation}
and
\begin{equation}C=xz+\frac{x^3z^2(t-1)(z+x(t-1)+x^2(t-1)z)}{1-xz(t-1)-x^2z^2(t-1)}+x^3(t-1)z.\end{equation}
\end{thm}

\subsection{The pattern \underline{132}}

Now we consider the pattern $\underline{132}.$
By Theorem \ref{corrispondenze}, an occurrence of this pattern in $\pi \in \mathbf I_n(3412)$ corresponds to six possible subwords in $\Psi(\pi),$
namely, $DUY$ and $HUY',$ where $Y$ and $Y'$
can be any letters in $\{U,D,H\}.$ The occurrences of such words  correspond to the occurrences of $DU$ and $HU$.

Also in this case we use the cluster method.
Here we have $S=\{HU,DU\}.$ Note that the only possible clusters formed by these two words are $HU$ and $DU$ themselves. The first of these two clusters
reduces to an up step and has depth $0,$ the second one reduces to a horizontal step and has depth $-1.$
Hence we have
\begin{equation}A_H(x,t,k)=x^2t,\end{equation}
\begin{equation}A_U(x,t,z)=A'_U(x,t,z)=x^2tz\end{equation}
and
\begin{equation}A_D(x,t,z)=A'_D(x,t,z)=A'_H(x,t,z)=0.\end{equation}

Theorem \ref{teoremacluster}  allows us to determine $F_{\underline{132}}(x,t+1,z).$
After the substitution $t\leftarrow t-1,$ we get
\begin{thm}\begin{equation} F_{\underline{132}}(x,t,z)=\frac{2}{1-xz+x^2t-x^2+\sqrt{(1-xz-x^2(t-1))^2-4x(x+x^2z(t-1))}}.\end{equation} \end{thm}

Notice that this generating function in the case $z=1$ encodes the distribution of \emph{weakly descending subpaths} over the set of Motzkin paths (see sequence A114690 in \cite{Sl}), where
a weakly descending subpath is a maximal subword consisting of $H$ and $D$ steps. In fact, every occurrence either of $HU$ or $DU$ breaks a weakly descending subpath. Hence, in every Motzkin path the number of weakly descending subpaths equals the number of occurrences of these two patterns increased by one.

\subsection{The pattern \underline{321}}

Theorem \ref{corrispondenze} shows that the occurrences of the  pattern $\underline{321}$ in $\pi$ correspond to the occurrences of $UUX,$
$X'DD $ and $UHD$ in $\Psi(\pi),$ where $X$ and $X'$ can be arbitrary letters in $\{U,H,D\}.$ Hence the occurrences of 
this pattern corresponds to the occurrences of $UU, DD$ and $UHD$ in $\Psi(\pi).$ 

Let $F_{\underline{321}}(x,t,z)$ be the corresponding generating function.

Denote by $X_w$ and $Y_w$ be the first and last step of a Motzkin path $w.$
We define  
\begin{itemize}
\item $A(x,t,z)$ to be the generating function of the set of  Motzkin paths such that $(X_w,Y_w)=(U,D)$,
\item $B(x,t,z)$ to be the generating function of the set of  Motzkin paths such that  $(X_w,Y_w)$ is either $(U,H)$ or $(H,D),$
\item $C(x,t,z)$ to be the generating function of the set of  Motzkin paths such that $(X_w,Y_w)=(H,H)$.
\end{itemize}
Observe that the first return decomposition implies that 
\begin{equation}\label{ricor}
F_{\underline{321}}=1+xzF_{\underline{321}}+x^2(1+xzt+Bt+At^2+C)F_{\underline{321}}.\end{equation}
A simple inclusion-exclusion argument yields
\begin{equation} A=F_{\underline{321}}-2xzF_{\underline{321}}+(xz)^2F_{\underline{321}}+xz-1.\end{equation}
Moreover it is easily seen that
\begin{equation} B=2A\frac{xz}{1-xz}=2(F_{\underline{321}}-2xzF_{\underline{321}}+(xz)^2F_{\underline{321}}+xz-1)\frac{xz}{1-xz}=2xz(F_{\underline{321}}-xzF_{\underline{321}}-1)\end{equation}
and \begin{equation} C=(xz)^2F_{\underline{321}}.\end{equation}

Substituting in (\ref{ricor}), we get  \begin{thm} $F_{\underline{321}}$  satisfies the following functional equation \begin{equation} aF_{\underline{321}}^2+bF_{\underline{321}}+c=0,\end{equation} where
\begin{equation} a= 2x^3zt-2x^4z^2t+x^2t^2-2x^3zt^2+x^4z^2t^2+x^4z^2,\end{equation}
\begin{equation} b=-1+xz-x^3zt-x^2t^2+x^2+x^3zt^2, \qquad c=1.\end{equation}\end{thm}

\subsection{The pattern \underline{312}}

This pattern correspond to occurrences of $UHH$ and $UHU$ in $\Psi(\pi).$
Let $F(x,t,z)$ be the corresponding generating function.  

Set  
\begin{equation} G(x,t_1,t_2,z)=\sum_n \sum_{d \in \mathcal M_n}x^n t_1^{o(UH)}t_2^{o(UHD)}z^{o(H)},\end{equation}
where $o(UH), o(UHD)$ and $o(H)$ denote the number of occurrences of the subwords $UH,UHD$ and $H$ in $d$.

By the first return decomposition we get the following recurrence for $G$.
\begin{equation}G=1+xzG+x^2G+x^3t_1t_2zG+x^3zt_1G(G-1)+x^2G(G-xzG-1).\end{equation}
In fact, a Motzkin path can either
\begin{itemize}
\item be empty, or
\item start by $UD,$ $UHD$, or 
\item be of the form $UHmDd$ or $Um'Dd$, where $m$ is a non empty Motzkin path, $m'$ is a non empty path starting with $U$, and $d$  an arbitrary path.
\end{itemize}
Hence we  get a functional equation satisfied by $F_{\underline{312}}(x,t,z)$ substituting in the previous equation $t_1\leftarrow t$ and $t_2\leftarrow \frac{1}{t}:$

\begin{equation} F_{\underline{312}}=1+xzF_{\underline{312}}+x^2F_{\underline{312}}+x^3zF_{\underline{312}}+x^3ztF_{\underline{312}}(F_{\underline{312}}-1)+x^2F_{\underline{312}}(F_{\underline{312}}-xzF_{\underline{312}}-1),\end{equation}
namely,
\begin{thm}The generating function $F_{\underline{312}}$ satisfies the following functional equation \begin{equation} aF_{\underline{312}}^2+bF_{\underline{312}}+c=0,\end{equation}
where
\begin{equation}a= x^3zt+x^2-x^3z,\end{equation}
\begin{equation}b=xz+x^2+x^3z-x^3zt-x^2-1, \qquad c=1.\end{equation}\end{thm}

\section{Inversions and descents over $\mathbf S_n(132,\underline{123})$}

We now turn to the case of permutations in $\mathbf S_n(132,\underline{123}).$

First of all we recall that, given a permutation $\pi\in\mathbf S_n(132,\underline{123}),$ if $\pi=w_1\ldots w_k$ is the decomposition of $\pi$ into ascending runs, then the $w_i's$ have length at most $2$ and the sequence of the heads of $\pi$ is a decreasing sequence. Moreover, the inverse of the map $\Gamma$ has an easy description in terms of tunnels of the Motzkin path, as in the case of the map $\Psi$. 

\begin{prop}\label{tunnel}
Let $d$ be a Motzkin path and $\pi$  the corresponding permutation in $\mathbf S_n(132,\underline{123}).$ 
Let $t_1t_2\ldots t_k$ be the sequence of tunnels of $d,$ listed in decreasing order of the $x$-coordinate of their leftmost point. The decomposition of $\pi$ into ascending runs is $\pi=w_1w_2\ldots w_k$ with $w_i=x_ix'_i,$ where $x_i$ is the $x$-coordinate of the first point of $t_i,$ increased by one, and $x'_i$ is the $x$-coordinate of the last point of $t_i.$  
\end{prop}

As an example consider the following Motzkin path
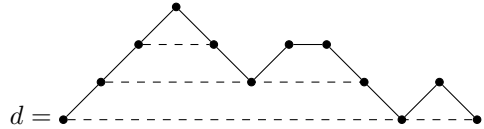
\begin{displaymath} d=
\begin{tikzpicture}%[every node/.style={draw,shape=circle,minimum size=1mm, inner sep=0mm, fill=black}]
\node[draw,shape=circle,minimum size=1mm, inner sep=0mm, fill=black] (11) at (11/2,0) {};
 \node[draw,shape=circle,minimum size=1mm, inner sep=0mm, fill=black] (10) at (10/2,1/2) {};
 \node[draw,shape=circle,minimum size=1mm, inner sep=0mm, fill=black] (9) at (9/2,0) {};
 \node[draw,shape=circle,minimum size=1mm, inner sep=0mm, fill=black] (8) at (8/2,1/2) {};
   \node[draw,shape=circle,minimum size=1mm, inner sep=0mm, fill=black] (7) at (7/2,1) {};
   \node[draw,shape=circle,minimum size=1mm, inner sep=0mm, fill=black] (6) at (6/2,1) {};
   \node[draw,shape=circle,minimum size=1mm, inner sep=0mm, fill=black] (5) at (5/2,1/2) {};
   \node[draw,shape=circle,minimum size=1mm, inner sep=0mm, fill=black] (4) at (4/2,1) {};
   \node[draw,shape=circle,minimum size=1mm, inner sep=0mm, fill=black] (3) at (3/2,3/2) {};
   \node[draw,shape=circle,minimum size=1mm, inner sep=0mm, fill=black] (2) at (2/2,1) {};
   \node[draw,shape=circle,minimum size=1mm, inner sep=0mm, fill=black] (1) at (1/2,1/2) {};
   \node[draw,shape=circle,minimum size=1mm, inner sep=0mm, fill=black] (0) at (0,0) {};
   \draw[] (0) -- (1) -- (2) -- (3) -- (4) -- (5) -- (6) -- (7) --(8)--(9)--(10)--(11);
    \draw[dashed] (0) -- (11)    (1) -- (8)    (2)--(4);
\end{tikzpicture}\end{displaymath}

The sequence of tunnels of $d$ is given by 9-11, 6-7, 5-8, 2-4, 1-5, 0-9, where each tunnel is represented by the $x$-coordinates of its first and last point.
Hence the corresponding permutation  is $\pi=10\,11\,7\,6\,8\,3\,4\,2\,5\,1\,9.$

Recall that a \emph{coinversion} in a permutation $\pi$ is a pair $(i,j)$ such that $i<j$ and $\pi_i<\pi_j.$ 
The number of coinversions of a permutation $\pi$ will be denoted by $\textrm{coinv}(\pi).$
Of course a permutation $\pi$ has $k$ coinversions if and only if it has $\binom{n}{2}-k$ inversions. 

Now we are interested in the generating function for permutations in $\mathbf S(132,\underline{123}):=\cup_{n\geq 0}  S_n(132,\underline{123})$
enumerated by number of coinversions and number of descents:
\begin{equation}  F(x,y,z)=\sum_{\pi \in \mathbf S(132,\underline{123})} x^ny^{\textrm{coinv}(\pi)}z^{\textrm{des}(\pi)}.\end{equation}
We have the following. 

\begin{prop}
Let $\pi\in \mathbf S_n(132,\underline{123})$ and let $\Gamma(\pi)$ be the corresponding Motzkin path. Then 
\begin{itemize}
\item $\textrm{coinv}(\pi)$ is the area of $\Gamma(\pi)$ and 
\item $\textrm{des}(\pi)$ is equal to the number of tunnels of $\Gamma(\pi)$ minus one. 
\end{itemize}
\end{prop}
\proof
Let $(i,j)$ be a coinversion of the permutation $\pi.$ Since the sequence of heads of $\pi$ is decreasing, $\pi_j$ is a tail, hence it corresponds to a down step in $\Gamma(\pi).$
Furthermore, given a down step $\bar D$ in $\Gamma(\pi)$ at position $k,$ consider the up step $\bar U$ that forms a tunnel with $\bar D,$ and denote by $h$ the position of $\bar U$. Then, by the construction of the map $\Gamma,$ the coinversions of $\pi$ having $k$ as  second element are precisely $(x,k)$ where $h\leq x\leq k-1.$ 
The number of such elements equals  the area of the trapezoid determined by the tunnel between $\bar U$ and $\bar D$. 

The second statement follows immediately form the fact that every descent in $\pi$ corresponds to a non initial head or head-tail.
\endproof

The above Proposition and the first return decomposition for Motzkin paths yield the following recurrence equation for the generating function $F.$

\begin{equation}\label{ricor2}
\begin{split}
F(x,y,z) & =1+x+zx(F(x,y,z)-1) +yz^2x^2(F(x,y,xz)-1)(F(x,y,z)-1)\\
&  +yzx^2(F(x,y,z)-1)+yzx^2(F(x,y,xz)-1) +yx^2. 
\end{split}
\end{equation}

Notice that $F(x,1,z)$ is the generating function of sequence A107131 in \cite{Sl}, while $F(x,y,1)$ is the generating function of sequence A129181 in \cite{Sl}.

\section{The distribution of consecutive patterns in $\mathbf S_n(132,\underline{123})$}

Now we enumerate permutations $\pi\in \mathbf S_n(132,\underline{123})$ by the number of occurrences of a consecutive pattern of length three. Needless to say, we  consider only the patterns 
$\underline{213},$ $\underline{231},$ $\underline{312}$ and $\underline{321}.$ 

\subsection{The pattern \underline{213}}

Let $F_{\underline{213}}(x,t)$ be the generating function of permutations $\pi\in\mathbf S_n(132,\underline{123})$ enumerated by length and number of occurrences of $\underline{213}$.

Note that an occurrence of this pattern in a permutation $\pi$ corresponds to
an occurrence in $\Gamma(\pi)$ of a sequence of the form $U\alpha D,$ where $\alpha$ is any non-empty Motzkin path.
We call such a sequence a \emph{long tunnel}.

In fact, an occurrence of $\underline{213}$ in $\pi$ is a sequence of consecutive letters $bac,$ with $a< b<c.$  Here, $ac$ is an ascending run $w_{i+1}$, while $b$ is either the tail or the head-tail of the preceding ascending run $w_i$.  

By Proposition \ref{tunnel}, $w_i$ and $w_{i+1}$ correspond to two tunnels $t_i,t_{i+1}$ such that $t_i$ lies above $t_{i+1}.$ Hence, the occurrence $bac$ of the pattern $\underline{213}$ corresponds to the long tunnel $t_{i+1}$.

Let $\widehat F(x,t,y)$ be the generating function for Motzkin paths enumerated by length ($x$),
occurrences of long tunnels ($t$) and peaks ($y$), \textit{i.e.}, occurrences of the sequence $UD.$

Notice that each non-empty Motzkin path can be either  a horizontal step followed by any Motzkin path,  or a peak followed by any Motzkin path, or a long tunnel followed by any Motzkin path. Hence, the  generating function $\widehat F(x,t,y)$ satisfies \begin{equation}\widehat F(x,t,y)=1+x\widehat F(x,t,y) + x^2y\widehat F(x,t,y) +x^2t(\widehat{F}(x,t,y)-1)\widehat{F}(x,t,y).\end{equation}

With the substitution $y \leftarrow 1$ we get

\begin{thm}The generating function $F_{\underline{213}}(x,t)$ satisfies the following equation: 
\begin{equation}aF_{\underline{213}}^2+bF_{\underline{213}}+1=0, \end{equation}
where 
\begin{equation}a=x^2t,\qquad \mbox{and }\qquad b=-1+x+x^2-x^2t.\end{equation}\end{thm}

\subsection{The pattern \underline{231}}

By Proposition \ref{tunnel} an occurrence  in $\pi$ of the pattern $\underline{231}$ corresponds to an occurrence in $\Gamma(\pi)$ of an up step in a non-initial position.
Let $\widehat F(x,t,y)$ be the generating function of Motzkin paths enumerated by length ($x$), number of non-initial up steps ($t$),
number of initial up steps ($y$).

We have the following recurrence for $\widehat F(x,t,y):$

\begin{equation}\label{ricorrenza1}
\widehat F(x,t,y)=1+x\widehat F(x,t,t)+x^2y\widehat F(x,t,t). 
\end{equation}

In fact, every non-empty Motzkin $d$ path can be decomposed either as 
$Hm,$ or $Um'Dm'',$
where $m,$ $m',$ and $m''$ are arbitrary Motzkin paths. Note that each up step  in $m,$ $m',$  or $m''$
cannot be at the initial position of $d.$

Substituting $y\leftarrow t$ in (\ref{ricorrenza1}) we find an expression for $\widehat F(x,t,t):$
\begin{equation}\widehat F(x,t,t)=\frac{1-x-\sqrt{1+x^2-2x-4x^2t}}{2x^2t}.\end{equation}
Substituting this expression in (\ref{ricorrenza1}) and then replacing  $y\leftarrow 1$ we get an expression for the generating function for permutations 
$\pi\in\mathbf S_n(132,\underline{123})$ enumerated by length ($x$) and number of occurrences of $\underline{231}$
($t$):
\begin{equation}F_{\underline{231}}(x,t)=\widehat F(x,t,1)=1+\widehat F(x,t,t) (x+x^2\widehat F(x,t,t)).\end{equation}

\subsection{The pattern \underline{312}}

An occurrence of the pattern $\underline{312}$ corresponds to an occurrence in $\Gamma(\pi)$ of a peak $p=UD$
such that $\Gamma(\pi)=\alpha p \beta$ where $\beta \neq D^k,$ $k\geq 0.$ We call such peak a \emph{non-final peak}.

Let $\widehat F(x,t,y)$ be the generating function for Motzkin paths enumerated by length ($x$), number of non-final peaks ($t$),
number of final peaks ($y$).

The first return decomposition implies that 

\begin{equation}
\begin{split}
 \widehat F(x,t,y)  & = \\
& =1+x\widehat F(x,t,y) +x^2y+x^2(\widehat F(x,t,y)-1)(\widehat F(x,t,t)-1) \\
& +x^2(\widehat F(x,t,y)-1)+x^2t(\widehat F(x,y,t)-1).  
\end{split} \end{equation}

Using the same arguments of the previous Subsection we get the following expression for the generating function $\widehat F(x,t,1):$ 
\begin{equation}F_{\underline{312}}(x,t)=\widehat F(x,t,1)=\frac{1-x^2\widehat F(x,t,t)+x^2-x^2t}{1-x-x^2\widehat F(x,t,t)-x^2t}\end{equation}
where \begin{equation}\widehat{F}(x,t,t)=\frac{-b-\sqrt{b^2-4a}}{2a}\end{equation}
with \begin{equation}a=x^2, \qquad \mbox{and}\qquad  b=-1+x-x^2+x^2t.\end{equation}

\subsection{The pattern \underline{321}}

An occurrence of the pattern $\underline{321}$ corresponds to an occurrence in $\Gamma(\pi)$ of a horizontal step
that is neither  in the first nor in  last position nor followed only by down steps. We call such a horizontal step a \emph{distinguished} horizontal step.

Let $\widehat F(x,t,y,z)$ be the generating function for Motzkin paths enumerated by length ($x$), number of distinguished  horizontal steps  ($t$),
number of horizontal steps in the first position ($y$), number of horizontal steps followed only by a (possibly empty) sequence of down steps ($z$). 

We have the following recurrence for $\widehat F(x,t,y,z):$

\begin{equation}
\begin{split}
\widehat F(x,t,y,z) & =\\
& = 1+xy(\widehat F(x,t,t,z)-1)+xz +x^2\widehat F(x,t,t,z) \\
& + x^2\widehat F(x,t,t,t)(\widehat F(x,t,t,z)-1)
\end{split}
\end{equation}

As a consequence
\begin{equation}F_{\underline{321}}(x,t)=\widehat F(x,t,1,1)=\frac{1-xt-x^2G+x-x^2t+x^2-x^3t+x^3}{-x^2+1-xt-x^2G}\end{equation}
where \begin{equation}G=\frac{-b-\sqrt{b^2-4a}}{2a}\end{equation}
with \begin{equation}a=x^2, \qquad \mbox{and} \qquad b=-1+xt.\end{equation}

\section{Appendix}\label{appendice}

In this appendix we describe the method that we used in Section \ref{involutionspatterns} to count Motzkin paths by 
occurrences of a set of given subpatterns. This method is a slight modification of the \emph{Goulden-Jackson cluster method}
used to enumerate words over an arbitrary finite alphabet  by occurrences of given subwords (\cite[p. 128]{Goulden}). In this context
the Goulden-Jackson cluster method does not apply directly, since the words  we are considering correspond to Motzkin paths, hence, they have particular constraints. 
Our method is inspired by those presented in \cite{wang}, where the author uses similar ideas to count Dyck words by occurrences of given subwords. 

Let $\mathcal A$ be the set of words in the alphabet $\{U,D,H\},$ and let $S\subseteq \mathcal A.$

Given $w\in \mathcal A,$ let $|w|$ be the length of $w,$
$|w|_L$  its number of steps of type  $L\in\{U,D,H\}$ and  $|w|_S$  the total number of occurrences in $w$
of subwords from $S.$ 

A \emph{marked subword} of a word $w=a_1\ldots a_n\in\mathcal A$ with respect to the set $S$ is a pair $(i,v)$
where $i$ is a positive integer, $v=a_ia_{i+1}\ldots a_{i+|v|-1}$ and $v\in S.$  A \emph{marked word} is a word $w\in\mathcal A$ with a (possibly empty)
set of marked subwords of $w.$ A \emph{cluster} with respect to $S$ is a marked word that is not the concatenation of two nonempty
marked words. 

As an example, consider  $S=\{UUU,DHU\}.$ The marked subwords of the 
word $w=UUUHUDDUUUDDDDHUHDD$ are $(1,UUU),$ $(8,UUU)$ and $(14,DHU).$
Hence \begin{displaymath}(w,\{(1,UUU),(14,DHU)\})\end{displaymath} is an example of a marked word.

Two clusters for $S$ are the marked words \begin{displaymath}(DHUUUUU,\{(1,DHU),(3,UUU),(4,UUU), (5,UUU)\})\end{displaymath}
\begin{displaymath}\begin{tikzpicture}[   txt/.style = {text height=2ex, text depth=0.25ex}]
  \node [txt]  (0) at (0,0) {D};
  \node [txt]  (0) at (1/2,0) {H};
  \node [txt]  (0) at (2/2,0) {U};
  \node [txt]  (0) at (3/2,0) {U};
  \node [txt]  (0) at (4/2,0) {U};
  \node [txt]  (0) at (5/2,0) {U};
  \node [txt]  (0) at (6/2,0) {U};
\draw [decorate,decoration={brace,mirror,amplitude=5pt}]
    (-0.2,-0.2) -- (1.2,-0.2) ;
    \draw [decorate,decoration={brace,mirror,amplitude=5pt}]
    (1.3,-0.2) -- (2.7,-0.2) ;
   \draw [decorate,decoration={brace,amplitude=5pt}]
    (0.8,0.5) -- (2.2,0.5) ;   
   \draw [decorate,decoration={brace,amplitude=5pt}]
    (1.8,0.3) -- (3.2,0.3) ;  
 \end{tikzpicture}\end{displaymath}
 and
\begin{displaymath}(DHUUUUU,\{(1,DHU),(3,UUU),(5,UUU)\})\end{displaymath}
\begin{displaymath}\begin{tikzpicture}[   txt/.style = {text height=2ex, text depth=0.25ex}]
  \node [txt]  (0) at (0,0) {D};
  \node [txt]  (0) at (1/2,0) {H};
  \node [txt]  (0) at (2/2,0) {U};
  \node [txt]  (0) at (3/2,0) {U};
  \node [txt]  (0) at (4/2,0) {U};
  \node [txt]  (0) at (5/2,0) {U};
  \node [txt]  (0) at (6/2,0) {U};
\draw [decorate,decoration={brace,mirror,amplitude=5pt}]
    (-0.2,-0.2) -- (1.2,-0.2) ;
   \draw [decorate,decoration={brace,amplitude=5pt}]
    (0.8,0.5) -- (2.2,0.5) ;   
   \draw [decorate,decoration={brace,amplitude=5pt}]
    (1.8,0.3) -- (3.2,0.3) ;  
 \end{tikzpicture}\end{displaymath}
whereas \begin{displaymath}(DHUUUUU,\{(1,DHU),(4,UUU),(5,UUU)\})\end{displaymath}
 \begin{displaymath}\begin{tikzpicture}[   txt/.style = {text height=2ex, text depth=0.25ex}]
  \node [txt]  (0) at (0,0) {D};
  \node [txt]  (0) at (1/2,0) {H};
  \node [txt]  (0) at (2/2,0) {U};
  \node [txt]  (0) at (3/2,0) {U};
  \node [txt]  (0) at (4/2,0) {U};
  \node [txt]  (0) at (5/2,0) {U};
  \node [txt]  (0) at (6/2,0) {U};
\draw [decorate,decoration={brace,mirror,amplitude=5pt}]
    (-0.2,-0.2) -- (1.2,-0.2) ;
    \draw [decorate,decoration={brace,mirror,amplitude=5pt}]
    (1.3,-0.2) -- (2.7,-0.2) ;
   \draw [decorate,decoration={brace,amplitude=5pt}]
    (1.8,0.3) -- (3.2,0.3) ;  
 \end{tikzpicture}\end{displaymath}
is not a cluster, because it can 
be seen as the juxtaposition of the marked words\\
 $(DHU,\{(1,DHU)\})$ and $(UUUU,\{(1,UUU),(2,UUU)\}).$

Note that, if $w'\in S,$  $(w',\{(1,w')\})$ is trivially a cluster.

We say that a word $w\in\mathcal A$ \emph{reduces to an up step} (\emph{to a down step}, \emph{to a horizontal
step})if $|w|_U-|w|_D=1$ ($-1,0,$ respectively). If one of these three cases occur we say that $w$ \emph{reduces to a single
step}. A cluster \emph{reduces to a single step}  if the underlying word does.

Now we  define the \emph{depth} of a word $w\in\mathcal A$ that reduces to a single step.
Draw the word $w$ in the lattice plane starting at the origin and assigning to the letters $U,D,H$ the usual steps.
Let $k$ be the  minimal $y$-coordinate reached by the resulting path.
We say that $w$ \emph{has depth} $k$ if it reduces to an up or horizontal step and $k+1$
if it reduces to a down step.
For example, the word $DDU$ has depth $-1,$ the word $DU$ has depth $-1,$ and the word $UDU$ has depth $0.$

Recall that the height of a step $d_i$ of a Motzkin path is 

\begin{displaymath}\begin{cases}
\mbox{the $y$-coordinate of the starting point of the step } d_i & \mbox{ if $d_i$ is either  $U$ or $H$}, \\
\mbox{the $y$-coordinate of the ending point of the step } d_i & \mbox{otherwise}.
\end{cases}\end{displaymath}
Observe that if a Motzkin path can be decomposed as $Um'Dm''$ then all the steps in $m'$  have height at least 1.

\begin{thm}\label{teoremacluster}
Let $S=\{w_1,\ldots,w_k\}$ be a subset of $\mathcal A$ such that  no words $w_i$ are proper subwords of other words in $S.$
Suppose that each  cluster formed by these words reduces to a single step and has depth greater than or equal to $-1.$ Let $A_H(x,t,z)$ be the generating function
of clusters that reduce to a horizontal step
enumerated by length ($x$), occurrences of $w_1,\ldots,w_k$ as subwords ($t$),
and horizontal steps ($z$). Denote by $A'_H(x,t,z)$
 the generating function of clusters that reduce to a horizontal step with depth $0.$
The generating functions $A_D(x,t,z),$ $A'_D(x,t,z),$ and $A_U(x,t,z),$ $A'_U(x,t,z)$ are defined in the same way for clusters
that reduce to a down step and an up step, respectively.
Then the generating function $F(x,t,z)$ for Motzkin paths enumerated by length, occurrences of the
words $w_1,\ldots,w_k$ and number of horizontal steps satisfies
\begin{equation}\label{equazgenfun}
F(x,t+1,z)= \frac{2ys}{2(1-l')ys-y's'+y's'l+y's'\sqrt{(1-l)^2-4ys}}, 
\end{equation}
where 
\begin{itemize}
\item $l=xz+A_H(x,t,z),$
\item $l'=xz+A'_H(x,t,z),$
\item $y=x+A_U(x,t,z),$
\item $y'=x+A'_U(x,t,z),$
\item $s=x+A_D(x,t,z),$
\item $s'=x+A'_D(x,t,z).$
\end{itemize}
\end{thm}
\proof
\begin{equation}F(x,t+1,z)=\sum_{w\in\mathcal M}x^{|w|}z^{|w|_H}(t+1)^{|w|_S}=\end{equation}
\begin{equation}\sum_{w\in\mathcal M}x^{|w|}z^{|w|_H}\sum_{k\geq 0}\binom{|w|_S}{k}t^k= \end{equation}
\begin{equation}\sum_{w\in\mathcal M}x^{|w|}z^{|w|_H}\sum_{T\subseteq S_w} t^{|T|},\end{equation}
where $S_w$ is the set of words in  $S$ contained in $w$ as subwords.

Hence  $F(x,t+1,z)$ counts Motzkin words
weighted by the number of marked subwords contained therein, by length and number of horizontal steps.

We want to show that  the right-hand side of (\ref{equazgenfun}) counts the same objects.

Let $G(y,s,l,y',s',l')$ be the generating function for Motzkin paths enumerated by occurrences
of $U,D$ and $H$ at non-zero height and by occurrences
of $U,D$ and $H$ at zero height. Hence, the formal power series $G_1(y,s,l):=G(y,s,l,y,s,l)$ is the generating
function of Motzkin paths enumerated by occurrences
of $U,D$ and $H.$

By the first return decomposition we get immediately  \begin{equation}G_1(y,s,l)=1+lG_1(y,s,l)+ysG_1^2(y,s,l)\end{equation}
and \begin{equation}G(y,s,l,y',s',l')=1+l'G(y,s,z,y',s',l')+y's' G(y,s,l,y',s',l')\cdot G_1(y,s,l).\end{equation}
As a consequence \begin{equation}G_1(y,s,l)=\frac{1-l-\sqrt{(1-l)^2-4ys}}{2ys}\end{equation}
and 
\begin{equation}
\begin{split}
G(y,s,l,y',s',l')&=\frac{1}{1-l'-y's'G_1(y,s,l)}\\
&=\frac{2ys}{2(1-l')ys-y's'+y's'l+y's'\sqrt{(1-l)^2-4ys}}.
\end{split}
\end{equation}

Let $\widehat G(x,t,z)$ be the generating function obtained from $G(y,s,l,y',s',l')$ by replacing
\begin{itemize}
\item the variable $l$ by $xz+A_H(x,t,z)$
\item the variable  $l'$ by $xz+A'_H(x,t,z)$
\item the variable  $y$ by $x+A_U(x,t,z)$
\item the variable  $y'$ by $x+A'_U(x,t,z)$
\item the variable  $s$ by $x+A_D(x,t,z)$
\item the variable  $s'$ by $x+A'_D(x,t,z).$
\end{itemize}

Note that $\widehat G(x,t,z)$ is precisely the right-hand side of Equation (\ref{equazgenfun}).

Let $w$ be a Motzkin word. Choose in $w$ some clusters $c_1,\ldots ,c_k.$ 
%each
%of the form $(w',\{(i_1,w_{j_1}),\ldots,(i_r,w_{j_r})),$ where $w'$ is a subword of $w$ and the $w_{j}$'s are elements of $S$. 
By hypothesis these clusters have depth $-1$ or $0.$

If in $w$ we replace each cluster  $c_i$  with the step that $c_i$ reduces to, we get another Motzkin word.

Conversely, given a Motzkin word $v$ we can choose in $v$ some up, down and horizontal steps, and replace them
by a cluster that reduces to an up, down and horizontal step, respectively, with the constraint that a step of height $0$ can be only replaced  by a cluster of depth $0.$

As a consequence the generating function $\widehat G(x,t,z)$
counts marked Motzkin words weighted by the number of marked subwords contained therein.

\endproof

\acknowledgements
We thank the anonymous referees for their detailed revisions and valuable suggestions.

\bibliographystyle{plain}
\bibliography{BIBLIOGRAFIA}

\end{document}